\documentclass[twocolumn,letterpaper, 10 pt,conference]{ieeeconf}

\IEEEoverridecommandlockouts

\usepackage{amssymb}
\usepackage{amsmath}
\usepackage{color}
\usepackage{graphics}
\usepackage{graphicx}
\usepackage{multirow}
\usepackage{cite}
\usepackage{algorithm}
\usepackage{savesym}
\savesymbol{AND}
\usepackage{epstopdf}
\usepackage{caption}
\usepackage{float}
\usepackage{algpseudocode}
\usepackage{times}
\usepackage{siunitx}

\usepackage{tikz}
\usetikzlibrary{arrows.meta}
\usetikzlibrary{shapes.misc}
\usetikzlibrary{shapes.arrows}
\usetikzlibrary{arrows}

\makeatletter
\def\set@curr@file#1{\def\@curr@file{#1}} 
\makeatother

\DeclareMathOperator{\spn}{span}

\newtheorem{theorem}{Theorem}
\newtheorem{remark}{Remark}
\newtheorem{lemma}{Lemma}
\newtheorem{corollary}{Corollary}
\newtheorem{definition}{Definition}
\newtheorem{proposition}{Proposition}

\title{\LARGE \bf Immersion-based model predictive control of constrained nonlinear systems: Polyflow approximation}

\author{Zheming Wang and Rapha\"el M. Jungers
\thanks{The authors are with the ICTEAM Institute, UCLouvain, Louvain-la-Neuve,1348, Belgium.}
\thanks{Rapha\"el M. Jungers is a  FNRS honorary Research Associate. This project has received funding from the European Research Council (ERC) under the European Union's Horizon 2020 research and innovation programme under grant agreement No 864017 - L2C. Rapha\"el M. Jungers is also supported by the Walloon Region and the Innoviris Foundation.}
\thanks{Email addresses:  zheming.wang@uclouvain.be (Zheming Wang), raphael.jungers@uclouvain.be (Rapha\"el M. Jungers)}}

\begin{document}
\maketitle

\begin{abstract}
In the framework of Model Predictive Control (MPC), the control input is typically computed by solving optimization problems repeatedly online. For general nonlinear systems, the online optimization problems are non-convex and computationally expensive or even intractable. In this paper, we propose to circumvent this issue by computing a high-dimensional linear embedding of discrete-time nonlinear systems. The computation relies on an algebraic condition related to the immersibility property of nonlinear systems and can be implemented offline. With the high-dimensional linear model, we then define and solve a convex online MPC problem. We also provide an interpretation of our approach under the Koopman operator framework. 
\end{abstract}


\section{Introduction}
Model Predictive Control (MPC) is a powerful control technology for general constrained systems due to its ability to handle hard constraints on control and states. The control action is determined online repeatedly by optimizing the future evolution of the trajectories of the system with the given cost function and constraints, see \cite{ART:MRRS00}. For general nonlinear systems, the online optimization problem to be solved is often non-convex and computationally expensive. To tackle this issue, many nonlinear MPC algorithms have been proposed (we refer the reader to the survey papers \cite{ART:C04,ART:M14}). In this paper, we intend to alleviate the online computation complexity by using linear embedding of nonlinear systems. More precisely, we use the state immersion method \cite{ART:MN83,ART:LM86,ART:LM88} to obtain a linear model with a higher dimension and convert the original non-convex MPC problem into a convex problem with this linear model. 


While the immersion method can be equivalent to feedback linearization in the special cases
where the immersion is a state diffeomorphism, they are in general different as an immersion does not necessarily preserve the dimension of the system. Hence, the immersion method has more freedom to find linear representations of nonlinear systems. It has been used in stabilization and observer design for nonlinear systems \cite{ART:KI83,ART:AO03}. Immersion is intrinsically related with the so-called “polynomial vector flows” or  \emph{polyflows} for continuous-time systems, see, e.g., \cite{ART:BM85,ART:VDE94,ART:VDE95}. Recently, a new immersion technique has been proposed in \cite{INP:JT19} for continuous-time autonomous systems using \emph{polyflows}. This technique is then extended to discrete-time autonomous systems in \cite{ART:WJ20}, which presents an immersibility algebraic condition that can be computed in a data-driven fashion. In this paper, we derive a similar immersibility algebraic  condition for control systems, which again allows to compute an immersion efficiently. 

Let us mention that linear embedding of nonlinear systems has already been used for analysis and control under the Koopman operator framework, see  \cite{ART:KKB17,ART:KM18,ART:MG20} and the references therein. The linear model is typically obtained by computing a truncation of the Koopman operator using a finite number of observables or basis functions. While our proposed \emph{polyflow} approximation can be viewed as a Koopman-based lifting approach with a special set of basis functions, it is derived from an immersibility algebraic condition and is equipped with an algorithmic way of determining the basis functions. See Section \ref{sec:comp} for more developments on the link between this work and  Koopman-based approaches.

The rest of the paper is organized as follows. The next section reviews some preliminary results on immersion and gives the problem statement. Section \ref{sec:optimal} presents an immersibility algebraic condition for control systems and discusses optimal control using linear embedding. In Section \ref{sec:comp}, we present a computational procedure for determining a linear model. Numerical experiments are provided in Section \ref{sec:Exp}.

\textbf{Notation}. The non-negative integer set is indicated by $\mathbb{Z}^+$.  $I_n$ is the $n\times n$ identity matrix and $\pmb{0}_{n\times m}$ is the $n\times m$ matrix of all zeros (subscript omitted when the dimension is clear). $\mathbb{B}_n$ is the unit closed ball in $\mathbb{R}^n$. $\|x\|_{p}$ denotes the $\ell_{p}$-norm of $x$($\|x\|=\|x\|_2$ by default). Given a map $T$, let $T(X)$ denote $\{T(x):x\in X\}$ and $T^{-1}(Y)$ denote the preimage of the set $Y$ under the map $T$, i.e., $T^{-1}(Y):=\{x:T(x)\in Y\}$ ($T$ is not necessarily invertible).  

\section{Preliminaries and problem statement}
We consider the following discrete-time nonlinear system
\begin{align}\label{eqn:fx}
x(t+1) = f(x(t),u(t)), ~~ t \in \mathbb{Z}^+
\end{align}
where $x(t)\in \mathbb{R}^n$ is the state vector, $u(t)\in \mathbb{R}^m$ is the input, and $f: \mathbb{R}^n \times \mathbb{R}^m  \rightarrow \mathbb{R}^n$ is some continuous function. Let $\pmb{u}$ denote a control sequence $\{u(0),u(1),\cdots,\}$ and $\phi(t;x,\pmb{u})$ denote the solution of (\ref{eqn:fx}) with a initial state $x$ at $t=0$ and a control sequence $\pmb{u}$.

\subsection{Immersion}
This paper is concerned with linear embedding of the nonlinear system (\ref{eqn:fx}). More precisely, we aim to immerse System (\ref{eqn:fx}) into a linear system ($\tilde{n} >n$):
\begin{align} \label{eqn:zAT}
\begin{aligned}
\Sigma(A,B, C) : \quad \tilde{x}(t+1) &= A\tilde{x}(t) + Bu(t),  \\
\tilde{y}(t) &= C\tilde{x}(t),  t \in \mathbb{Z}^+
\end{aligned}
\end{align}
where $\tilde{x}(t)\in \mathbb{R}^{\tilde{n}}$, $\tilde{y}(t)\in \mathbb{R}^n$, $A\in \mathbb{R}^{\tilde{n}\times \tilde{n}}$, and $C \in \mathbb{R}^{n\times \tilde{n}}$. Similarly, let $\tilde{\phi}(t;\tilde{x},\pmb{u})$ denote the solution of (\ref{eqn:zAT}) with a initial state $\tilde{x}$ (at time $t=0$) and a control sequence $\pmb{u}$.

Let us first recall the notion of immersibility, see, e.g., \cite{ART:MN83,ART:LM86,ART:LM88}.
\begin{definition}\label{def:immersion}
System (\ref{eqn:fx}) is immersible into  $\Sigma(A,B, C)$ given in (\ref{eqn:zAT}) if there exists a map $T:\mathbb{R}^n \rightarrow \mathbb{R}^{\tilde{n}}$ such that, for any initial state $x\in \mathbb{R}^n$ and control sequence $\pmb{u}$,  $\phi(t;x,\pmb{u}) = C\tilde{\phi}(t;T(x),\pmb{u})$ for all $t\in \mathbb{Z}^+$, where $\phi(t;x,\pmb{u})$ is the solution of (\ref{eqn:fx}) and $\tilde{\phi}(t;T(x),\pmb{u})$ is the solution of (\ref{eqn:zAT}). 
\end{definition}

From this definition, it can be verified that System (\ref{eqn:fx}) is immersible into  $\Sigma(A,B, C)$ if there exist a map $T: \mathbb{R}^n \rightarrow \mathbb{R}^{\tilde{n}}$ and matrices $A\in \mathbb{R}^{\tilde{n}\times \tilde{n}}, B\in \mathbb{R}^{\tilde{n}\times m}$ and $C\in \mathbb{R}^{n\times \tilde{n}}$ such that, $\forall (x,u) \in \mathbb{R}^n \times \mathbb{R}^m$,
\begin{align}\label{eqn:TfxuABC}
T(f(x,u)) =  AT(x)+Bu, ~~x = CT(x).
\end{align}
A necessary and sufficient condition of immersibility is given in \cite{ART:LM88}. To state this immersibility condition, we bring in the notation in \cite{ART:LM88}. Let 
\begin{align}\label{eqn:hatf}
\hat{f}^{\ell+1}(x,u) = f(\hat{f}^{\ell}(x,u),0), \ell \ge 1.
\end{align}
with $\hat{f}^1(x,u) =f(x,u)$ and $\hat{f}^0(x,u) =x$.  Define the space of $\mathbb{R}$-valued functions
\begin{align}\label{eqn:Phi}
\Phi = \spn\limits_{\mathbb{R}}\{&x_1,\cdots, x_n, \hat{f}^1_1(x,0),\cdots, \hat{f}^1_n(x,0),\cdots,\nonumber\\
&\hat{f}^k_1(x,0),\cdots, \hat{f}^k_n(x,0),\cdots \},
\end{align}
where $\hat{f}^k_i(x,0)$ is the $i^{th}$ component of $\hat{f}^k_n(x,0)$ for $i=1,2,\cdots,n$. Then, we have the following characterization of immersible systems:
\begin{theorem}[\cite{ART:LM88}]\label{thm:imm}
System (\ref{eqn:fx}) is immersible into a linear system in form of (\ref{eqn:zAT}) if and only if: i) $\Phi$, defined as in (\ref{eqn:Phi}), is finite-dimensional; 2) $D_{u^T} \hat{f}^\ell(x,u)$ is constant  (independent of $x$ and $u$) for all $\ell \ge 1$, where $\hat{f}^\ell(x,u)$ is defined as in (\ref{eqn:hatf}).
\end{theorem}

\subsection{Problem statement}
The goal of this paper is to design controllers using the linear embedding in (\ref{eqn:zAT}). We consider constrained nonlinear systems with state and input constraints as follows:
\begin{align}\label{eqn:xXuU}
x(t) \in X, ~~ u(t) \in U, ~~ t \in \mathbb{Z}^+
\end{align}
where $X\subseteq \mathbb{R}^n$ and $U\subseteq \mathbb{R}^m$ are the given state and input constraint sets. The control objective is to steer the state to the origin while fulfilling the constraints in (\ref{eqn:xXuU}). This can be accomplished under the MPC framework. Given the current state $x(t)$ at time $t$, we define the MPC problem as follow:
\begin{subequations}\label{eqn:Pxt}
\begin{align}
 \min_{\pmb{u}_N} &\sum\limits_{i=0}^{N-1} l(x_i,u_i) + l_f(x_N)\\
\textrm{s.t.   } \quad & x_{i+1} = f(x_i,u_i), x_i \in X, u_i \in U,\\
& x_0 = x(t), x_N \in X_f, i = 0,1,\cdots, N-1,
\end{align}
\end{subequations}
where $N$ is the horizon length, $\pmb{u}_N := \{u_0,u_1,\cdots, u_{N-1}\}$, $l(x,u)$ is the stage cost, $l_f(x)$ is the terminal cost, and $X_f \subseteq X$ is some appropriate terminal set. Once the problem above is solved, the first element of the control sequence is applied to the system and this process is repeated. Under some condition, recursive feasibility and stability is guaranteed, see, e.g., \cite{ART:MRRS00}. While there are well-known nonlinear MPC techniques, see \cite{ART:C04,ART:M14}, solving the nonlinear MPC problem online is computationally expensive in general. In this paper, we propose to use a lifted linear model (which is computed offline) in the online MPC problem to reduce computational burden, motivated by the immersibility property. Eventually, we migrate the complexity of the online computation of nonlinear optimization problems to the offline computation of linear embedding.



\section{Optimal control via immersion}\label{sec:optimal}
In this section, we discuss optimal control of nonlinear systems when a linear immersion is available.  

\subsection{Discrete-time polyflows}
We first show a nilpotency property of immersible discrete-time nonlinear systems.  For  continuous-time autonomous systems, when a system is immersible to a linear system, it enjoys a Lie derivative nilpotency property and the solution of the system is called a \emph{polyflow}, see \cite{ART:VDE94}. Based on this nilpotency property, a linearization technique called \emph{polyflow} approximation has been proposed in \cite{INP:JT19} and it is proved that this approximation exhibits better convergence properties than Taylor approximation. While \emph{polyflows} are defined for continuous-time systems, the initial definition of  \cite{INP:JT19} can be generalized to discrete-time systems, see \cite{ART:WJ20}. In this paper, we further generalize the concept of \emph{polyflows}  for control nonlinear systems in the form of (\ref{eqn:fx}).


Similar to the Lie derivative nilpotency property for continuous-time, we also derive a nilpotency condition for immersibility in discrete-time, which allows to use powerful linear optimal control techniques \cite{BOO:LVS12} for nonlinear systems. Before we present the nilpotency condition, the following lemma is needed.

\begin{lemma}\label{lem:fk0}
Consider System (\ref{eqn:fx}), let $\hat{f}^{\ell}(x,u)$ be defined as in  (\ref{eqn:hatf}) for all $\ell \in \mathbb{Z}^+$. Then, it holds that $\hat{f}^{k}(f(x,u),0) = \hat{f}^{k+1}(x,u)$ for all $k \in \mathbb{Z}^+$.
\end{lemma}
Proof: The proof goes by induction. Suppose $\hat{f}^{k}(f(x,u),0) = \hat{f}^{k+1}(x,u)$ holds for some integer $k \ge 0$. Then, $\hat{f}^{k+1}(f(x,u),0) = f(\hat{f}^{k}(f(x,u),0),0) = f(\hat{f}^{k+1}(x,u),0) =  \hat{f}^{k+2}(x,u)$, from the definition in (\ref{eqn:hatf}). It is obvious that $\hat{f}^{0}(f(x,u),0) = f(x,u) = \hat{f}^{1}(x,u)$.  $\Box$

Inspired the Lie derivative nilpotency property in \cite{ART:VDE94}, we modify Theorem \ref{thm:imm} and obtain a nilpotency condition for immersibility of discrete-time nonlinear systems.

\begin{theorem}\label{thm:immT}
Consider System (\ref{eqn:fx}), let $\hat{f}^{\ell}(x,u)$ be defined as in  (\ref{eqn:hatf}) for all $\ell \in \mathbb{Z}^+$. System (\ref{eqn:fx}) is immersible into a linear system in form of (\ref{eqn:zAT}) if and only if there exist a finite $k \in \mathbb{Z}^+$, $\{\alpha_\ell \in \mathbb{R}^{n\times n}\}_{\ell=0}^k$ and $\{b_\ell\in \mathbb{R}^{n\times m}\}_{\ell=0}^k$ such that
\begin{align}
\hat{f}^{k+1}(x,0) &= \sum_{\ell=0}^k \alpha_\ell \hat{f}^\ell(x,0),  \label{eqn:hatfkplu1x0}\\
\hat{f}^{\ell}(x,u) &=  \hat{f}^{\ell}(x,0) + b_{\ell} u, \ell = 1,2,\cdots, k+1. \label{eqn:hatfxu}
\end{align}
\end{theorem}
Proof: (Necessity) Suppose System (\ref{eqn:fx}) is immersible into a linear system in the form of (\ref{eqn:zAT}). From Theorem \ref{thm:imm}, $\Phi$ is finite-dimensional and there exists a $k\in \mathbb{Z}^+$ such that  $\{x_1,\cdots,  x_n, \hat{f}^1_1(x,0),\cdots, \hat{f}^1_n(x,0),\cdots,\hat{f}^k_1(x,0),\cdots,$ $\hat{f}^k_n(x,0)\}$ spans $\Phi$. Hence, there exist $\{\alpha_\ell \in \mathbb{R}^{n\times n}\}_{i=0}^k$ such that $\hat{f}^{k+1}(x,0) = \sum_{i=0}^k \alpha_\ell \hat{f}^\ell(x,0)$. Again, from Theorem \ref{thm:imm}, $D_{u^T} \hat{f}^\ell(x,u)$ is constant for all $\ell \ge 1$. Let $b^{\ell} = D_{u^T} \hat{f}^\ell(x,u)$ for $\ell \ge 1$.  Thus, $\hat{f}^\ell(x,u)$ can be written as $\hat{f}^\ell(x,u) = \hat{f}^\ell(x,0) + b^{\ell} u$ for all $\ell \ge 1$.\\
(Sufficiency) Suppose (\ref{eqn:hatfkplu1x0}) -- (\ref{eqn:hatfxu}) hold. From Definition \ref{def:immersion}, to show that System (\ref{eqn:fx}) is immersible into a linear system, we only need to show that there exist a map $T: \mathbb{R}^n \rightarrow \mathbb{R}^{\tilde{n}}$ and matrices $A\in \mathbb{R}^{\tilde{n}\times \tilde{n}}, B\in \mathbb{R}^{\tilde{n}\times m}$ and $C\in \mathbb{R}^{n\times \tilde{n}}$ such that (\ref{eqn:TfxuABC}) hold. Let 
\begin{align}
T(x) := \left( \begin{array}{c}
x\\
\hat{f}^1(x,0)\\
\vdots\\
\hat{f}^k(x,0)
\end{array} \right).
\end{align}
From (\ref{eqn:hatfkplu1x0}), we have
\begin{align}
&\left( \begin{array}{c}
\hat{f}^1(x,0)\\
\hat{f}^2(x,0)\\
\vdots\\
\hat{f}^{k+1}(x,0)
\end{array} \right) \nonumber\\
= &\left(\begin{array}{ccccc}
\pmb{0} & I_n & \pmb{0} & \cdots & \pmb{0}\\
\pmb{0} & \pmb{0} & I_n & \cdots & \pmb{0}\\
\vdots & \vdots & \vdots & \vdots & \vdots\\
\pmb{0} & \pmb{0} & \cdots & \pmb{0}  & I_n\\
\alpha_0 & \alpha_1 & \cdots & \alpha_{k-1} & \alpha_{k}
\end{array}
\right) \left( \begin{array}{c}
x\\
\hat{f}^1(x,0)\\
\vdots\\
\hat{f}^k(x,0)
\end{array} \right) \nonumber\\
 := & AT(x).
\end{align}
With this and Lemma \ref{lem:fk0}, we can get
\begin{align}
T(f(x,u)) &= \left( \begin{array}{c}
f(x,u)\\
\hat{f}^1(f(x,u),0)\\
\vdots\\
\hat{f}^k(f(x,u),0)
\end{array} \right) = \left( \begin{array}{c}
\hat{f}^1(x,u)\\
\hat{f}^2(x,u)\\
\vdots\\
\hat{f}^{k+1}(x,u)
\end{array} \right) \nonumber\\
& \stackrel{\text{(\ref{eqn:hatfxu})}}{=}
 \left( \begin{array}{c}
\hat{f}^1(x,0) + b^1 u\\
\hat{f}^2(x,0) + b^2 u\\
\vdots\\
\hat{f}^{k+1}(x,0) + b^{k+1} u
\end{array} \right) \nonumber \\
&= AT(x) +  \left( \begin{array}{c}
b^1 \\
b^2 \\
\vdots\\
b^{k+1}
\end{array} \right) u \nonumber\\
 :&= AT(x) +  B u
\end{align}
Therefore, (\ref{eqn:TfxuABC}) is satisfied with $C = [I_n 0_{n\times kn}]$. $\Box$ 

The nilpotency condition above is in a similar form of the Lie derivative nilpotency property for \emph{polyflows} in \cite{ART:VDE94}. Hence, analogously, we say that the solution of System (\ref{eqn:fx}) satisfying (\ref{eqn:hatfkplu1x0}) -- (\ref{eqn:hatfxu}) is a \emph{polyflow} when the condition in Theorem \ref{thm:immT} is satisfied.  
 
To  show the equivalence to Theorem \ref{thm:imm}, we derive the following corollary.
\begin{corollary}
Consider System (\ref{eqn:fx}), let $\hat{f}^{\ell}(x,u)$ be defined as in  (\ref{eqn:hatf}) for all $\ell \in \mathbb{Z}^+$. If there exist a finite $k\in \mathbb{Z}^+$, $\{\alpha_\ell \in \mathbb{R}^{n\times n}\}_{\ell=0}^k$ and $\{b_\ell\in \mathbb{R}^{n\times m}\}_{\ell=0}^{k}$ such that (\ref{eqn:hatfkplu1x0}) - (\ref{eqn:hatfxu}) hold, then, $D_{u^T} \hat{f}^\ell(x,u)$ is constant for all $\ell \ge 1$.
\end{corollary}
Proof: With (\ref{eqn:hatfkplu1x0}) -- (\ref{eqn:hatfxu}), we only need to show that $D_{u^T} \hat{f}^{k+2}(x,u)$ is constant. From Lemma \ref{lem:fk0}, we have $\hat{f}^{k+2}(x,u) = \hat{f}^{k+1}(f(x,u),0) = \sum_{i=0}^k \alpha_\ell \hat{f}^i(f(x,u),0) =\sum_{\ell=0}^k \alpha_\ell \hat{f}^{\ell+1}(x,u) =  \sum_{\ell=0}^k \alpha_\ell (\hat{f}^{\ell+1}(x,0) + b_{\ell+1} u)$. Hence, $D_{u^T} \hat{f}^{k+2}(x,u)$ is constant. $\Box$

While Theorem \ref{thm:immT} is a modification of Theorem \ref{thm:imm} (originated from \cite{ART:LM88}), the reformulation in (\ref{eqn:hatfkplu1x0})-(\ref{eqn:hatfxu}) facilitates the computation of linear embedding approximations, as shown in the next section.

\textbf{A special class of nonlinear systems}  \hspace{1 mm}  The immersibility property only hold for very special classes of systems. In a recent paper\cite{ART:WJ20}, we have shown a special class of discrete-time autonomous nonlinear systems with guaranteed immersibility. With linear control input, such systems become:
\begin{align}\label{eqn:x1x2}
x^1(t+1)  &= A^1 x^1(t) + \varphi(x^2(t)) + B^1 u(t), \nonumber\\
x^2(t+1) &= A^2 x^2(t), t \in \mathbb{Z}^+
\end{align}
where $x^1\in \mathbb{R}^{n_1},x^2\in \mathbb{R}^{n_2},n_1+n_2 = n, A^1\in \mathbb{R}^{n_1\times n_1},  A^2\in \mathbb{R}^{n_2\times n_2},B^1\in \mathbb{R}^{n_1\times m} $ and $\varphi: \mathbb{R}^{n_2} \rightarrow \mathbb{R}^{n_1}$ is some polynomial function. 

\begin{corollary}
System (\ref{eqn:x1x2}) is immersible to a linear system in form of (\ref{eqn:zAT}).
\end{corollary}
Proof: From Theorem 2 in \cite{ART:WJ20}, there exist a finite $k\in \mathbb{Z}^+$ and $\{\alpha_\ell \in \mathbb{R}^{n\times n}\}_{\ell=0}^k$ such that (\ref{eqn:hatfkplu1x0}) is satisfied. The condition (\ref{eqn:hatfxu}) can be easily verified. $\Box$

\subsection{Linear quadratic regulator with immersion}
Now, we consider optimal control of System (\ref{eqn:fx}) with no constraint by the use of immersion. Suppose System (\ref{eqn:fx}) is immersible into $\Sigma(A,B, C)$ under the map $T(x)$, we formulate the following regulation problem of System (\ref{eqn:fx}) with the initial state $x(0)$
\begin{subequations}
\begin{align}
&\min_{\pmb{u}} \sum_{t=0}^{\infty} \|x(t)\|_Q^2 + \|u(t)\|_R^2\\
\textrm{s.t} \quad &x(t+1) = f(x(t),u(t)), t \in \mathbb{Z}^+
\end{align}
\end{subequations}
where $Q \succeq 0$ and $R\succ 0$ are user-defined parameters. With the condition (\ref{eqn:TfxuABC}), this optimal control problem can be equivalently written as the linear–quadratic regulator (LQR) for $\Sigma(A,B, C)$:
\begin{subequations}
\begin{align}
&\min_{\pmb{u}} \sum_{t=0}^{\infty} \|\tilde{x}(t)\|_{C^TQC}^2 + \|u(t)\|_R^2\\
\textrm{s.t.} \quad & \tilde{x}(t+1) = A\tilde{x}(t)+Bu(t), \tilde{x}(0) = T(x(0)).
\end{align}
\end{subequations}
Hence, when $(A,B)$ is stabilizable and $(A,C)$ is observable, the optimal feedback control is 
\begin{align}
u(t) = KT(x(t)), t \in \mathbb{Z}^+, 
\end{align}
where $K$ is obtained by solving the discrete-time algebraic Riccati equation \cite{BOO:LVS12}
\begin{align}
P =& A^TPA - (A^TPB)(R+B^TPB)^{-1}(B^TPA) \nonumber\\
&+ C^TQC, \quad K = -(R+B^TPB)^{-1}B^TPA. \label{eqn:PK}
\end{align}

\subsection{Immersion-based MPC}
In the presence of constraints in (\ref{eqn:xXuU}), we can design MPC for System (\ref{eqn:fx}) using $\Sigma(A,B, C)$ and $T(x)$ in a similar way. We consider the following the stage and terminal costs in (\ref{eqn:Pxt}) 
\begin{align}
l(x,u) = \|x_i\|_Q^2 + \|u_i\|_R^2, l_f(x) = \|T(x)\|_P^2.
\end{align}
Let 
\begin{align}
\tilde{X}_f  := \{&x\in \mathbb{R}^{\tilde{n}}: C(A+BK)^kx \in X, \nonumber \\
&K(A+BK)^k x \in U, k\ge 0\}.
\end{align}
The terminal constraint set in (\ref{eqn:Pxt}) is chosen to be the preimage of $\tilde{X}_f $ under the map $T(x)$:
\begin{align}
X_f = T^{-1}(\tilde{X}_f).
\end{align}
The set $\tilde{X}_f $ is the maximal invariant set of the system $\tilde{x}^+=(A+BK)\tilde{x}$ and can be efficiently computed using the classic algorithm in \cite{ART:GT91} when $X$ and $U$ are polytopes. There also exist algorithms for handling non-convex constraints, see \cite{INP:AJ16,INP:WJO19,ART:WJO21Com}. As shown in \cite{ART:WJO20}, from the immersibility property, $T^{-1}(\tilde{X}_f)$ is the maximal invariant set of the closed system $x^+ = f(x,KT(x))$ with the control law $u=KT(x)$. With these definitions, the MPC problem is cast as follows:
\begin{subequations}\label{eqn:MPCT}
\begin{align}
\min_{\pmb{u}_N} &\sum\limits_{i=0}^{N-1} (\|x_i\|_Q^2 + \|u_i\|_R^2) + \|T(x_N)\|_P^2\\
\textrm{s.t.   } \quad & x_{i+1} = f(x_i,u_i), x_i \in X, u_i \in U,\\
&x_0 = x(t), x_N \in T^{-1}(\tilde{X}_f), i = 0,\cdots, N-1.
\end{align}
\end{subequations}
With (\ref{eqn:TfxuABC}), Problem (\ref{eqn:MPCT}) can be equivalently written as,
\begin{subequations}\label{eqn:MPCABC}
\begin{align}
\quad \min_{\pmb{u}_N} &\sum\limits_{i=0}^{N-1} (\|\tilde{x}_i\|_{C^TQC}^2 + \|u_i\|_R^2) + \|\tilde{x}_N\|_P^2\\
\textrm{s.t.   } \quad & \tilde{x}_{i+1} = A\tilde{x}_{i} + Bu_i, C\tilde{x}_i \in X, u_i \in U, \label{eqn:tildexiplus1AB}\\
&\tilde{x}_0 = T(x(t)), \tilde{x}_N \in \tilde{X}_f, i = 0,\cdots, N-1. \label{eqn:tildexTxtXf}
\end{align}
\end{subequations}
Let the problem above be denoted as $\tilde{\mathbb{P}}_{\Sigma(A,B,C)}(x)$. Under the assumption that $X$ and $U$ are convex, this problem is convex and can be efficiently solved. In particular, when $X$ and $U$ are convex polytopes, it becomes a quadratic optimization problem. The feasible domain of $\tilde{\mathbb{P}}_{\Sigma(A,B,C)}(x)$ is given by
\begin{align}
\mathcal{D}_{\Sigma(A,B,C)}:= \{x\in X: \exists \pmb{u}_N \textrm{ s.t. } (\ref{eqn:tildexiplus1AB})-(\ref{eqn:tildexTxtXf})\}.
\end{align}

\section{Linear embedding computation}\label{sec:comp}
This section presents the detailed procedure for computing linear embedding $\Sigma(A,B, C)$. We also show the connections between the proposed immersion-based approach  and  Koopman-based optimal control \cite{ART:KKB17,ART:KM18}.
\subsection{Polyflow approximation} 
For general nonlinear systems, the immersibility property does not hold and we can only obtain approximate immersion or linear embedding. 
Similar to the \emph{polyflow} approximation method in \cite{INP:JT19} for continuous-time autonomous systems, we develop a linearization technique for discrete-time control systems by using the nilpotency condition in Theorem \ref{thm:immT}. Given the similarity in the form of the linearized system, we also call it \emph{polyflow} approximation.

\begin{definition}
Given a $k \in \mathbb{Z}^+$ and $\{\alpha_\ell \in \mathbb{R}^{n\times n}\}_{\ell=0}^k$, let
\begin{align}\label{eqn:ABCpoly}
A &= \left(\begin{array}{ccccc}
\pmb{0} & I_n & \pmb{0} & \cdots & \pmb{0}\\
\pmb{0} & \pmb{0} & I_n & \cdots & \pmb{0}\\
\vdots & \vdots & \vdots & \vdots & \vdots\\
\pmb{0} & \pmb{0} & \cdots & \pmb{0}  & I_n\\
\alpha_0 & \alpha_1 & \cdots & \alpha_{k-1} & \alpha_{k}
\end{array}
\right), \nonumber\\
B &= \left( \begin{array}{c}
D_{u^T} \hat{f}^1(0,0)\\
D_{u^T} \hat{f}^2(0,0)\\
\vdots\\
D_{u^T} \hat{f}^{k+1}(0,0)
\end{array} \right), C = \left( \begin{array}{cc} I &
0_{n\times kn}
\end{array} \right).
\end{align}
We call $\Sigma(A,B, C)$ where matrices $(A,B,C)$ are in the form of (\ref{eqn:ABCpoly})   a $k^{th}$ \emph{polyflow} approximation of System (\ref{eqn:fx}).
\end{definition}

To determine a \emph{polyflow} approximation of order $k$, we need to find $\{\alpha_\ell \in \mathbb{R}^{n\times n}\}_{\ell=0}^k$ such that $\|\hat{f}^{k+1}(x,0) - \sum_{\ell=0}^k \alpha_\ell \hat{f}^\ell(x,0)\|$ is minimized for any $x\in X$. To do so, we generate a set of $M$ data points and solve the following problem for numerical \emph{polyflow} approximation
\begin{align}\label{eqn:alphakplus1}
\min_{\pmb{\alpha}_k} \sum_{i=1}^M \|\hat{f}^{k+1}(x_i,0) - \pmb{\alpha}_k\hat{F}^k(x_i)\|_2^2,
\end{align}
where $\pmb{\alpha}_k:=\{\alpha_0,\alpha_1,\cdots, \alpha_k\}$ and 
\begin{align}\label{eqn:Fkfhat}
\hat{F}^k(x) = \left( \begin{array}{c}
x\\
\hat{f}^1(x,0)\\
\vdots\\
\hat{f}^k(x,0)
\end{array} \right), k \in \mathbb{Z}^+.
\end{align}


\subsection{EDMD and polyflow basis}
Besides \emph{polyflow} approximation above, the nilpotency condition in Theorem \ref{thm:immT} also allows us to compute linear embedding in a fully data-driven fashion using the extended dynamical mode decomposition (EDMD) algorithm \cite{ART:WKR15}, like Koopman-based optimal control \cite{ART:KKB17,ART:KM18}.  We first need the following proposition.

\begin{proposition}\label{prop:TfABC}
System (\ref{eqn:fx}) is immersible into  $\Sigma(A,B, C)$, if and only if there exist a map $T: \mathbb{R}^n \rightarrow \mathbb{R}^{\tilde{n}}$ and matrices $A\in \mathbb{R}^{\tilde{n}\times \tilde{n}}, B\in \mathbb{R}^{\tilde{n}\times m}$ and $C\in \mathbb{R}^{n\times \tilde{n}}$ such that (\ref{eqn:TfxuABC}) holds  $\forall (x,u) \in \mathbb{R}^n \times \mathbb{R}^m$.
\end{proposition}

Although this result can be easily derived from Theorem \ref{thm:immT}, it can be proved independently. An elementary proof is given in the appendix. 

Motived by Proposition \ref{prop:TfABC}, to compute a linear approximation, we first select the map $T(x)$. Then, we generate a set of $M$ snapshot data points $\{(x_i,u_i,x_i^+):i=1,2,\cdots, M\}$( where $x_i^+ = f(x_i,u_i)$) from $X\times U$ and formulate the following problem, motivated from (\ref{eqn:TfxuABC}),
\begin{align}
\min_{A,B} \sum_{i=1}^M \|T(x_i^+) - AT(x_i)-Bu_i\|_2^2. \label{eqn:ABM}
\end{align}
Note that $C$ is usually known from the selection of $T(x)$. We arrive at the same problem as in \cite{ART:KKB17,ART:KM18}, but from a different path. This shows the link between the immersibility property and the Koopman operator, as already mentioned in \cite{ART:WJ20} for discrete-time autonomous systems.

In particular, when the dynamics $f(x,u)$ is known, we propose to use a special set of basis functions. We call the functions $\hat{F}^k(x)$ defined in (\ref{eqn:Fkfhat}) the \emph{polyflow} basis of degree $k$ for any $k \in \mathbb{Z}^+$. The \emph{polyflow} basis is also called delay coordinates for autonomous systems under the Koopman operator framework \cite{INP:SM15, ART:AM17}. To determine an appropriate basis, we solve Problem (\ref{eqn:alphakplus1}) and increase $k$ until the approximation error is sufficiently small. In case that there is redundancy in $\hat{F}^k(x)$, i.e., some component is a linear combination of others, we remove the redundant components and let the lifting map $T(x)$ be $V\hat{F}^k(x)$, where $V\in \mathbb{R}^{\tilde{n}\times n(k+1)}$ is some proper matrix and $\tilde{n}$ is the number of linearly independent components of $\hat{F}^k(x)$. The issue of redundancy occurs when one component of $f(x,0)$ is linear in $x$. In a data-driven fashion, we can also remove redundancy using singular value decomposition (SVD). When the system is equipped with a well-designed nominal controller $\kappa_{nominal}(x)$, we can also choose the basis from the closed-loop system $x^+=f(x,\kappa_{nominal}(x))$.

\begin{remark}
While the link between the immersibility property and the Koopman operator is mentioned above, 
we emphasize that our immersion-based approach has two key differences compared with Koopman-based approaches \cite{ART:KKB17,ART:KM18}. First, the underlying problem in (\ref{eqn:ABM}) stems from the immersibility condition in Theorem \ref{thm:immT} while Koopman-based approaches essentially solve a model reduction problem of the infinite-dimensional Koopman operator. Second, we provide a way to compute endogenously a basis of functions motivated by the nilpotency immersibility condition. On the contrary, the Koopman-based approaches in \cite{ART:KKB17,ART:KM18} fix the basis a priori, which incurs additional conservativeness.
\end{remark}

\section{Experiment: optimal pest control}\label{sec:Exp}
Consider the pest control problem in \cite{ART:DHCLG14}: 
\begin{align*}
v(t+1) &= v(t) + cv(t)(1-v(t)/\kappa)-rv(t)p(t), \\
p(t+1) &= dp(t)+ v(t)p(t)-a(t)p(t),
\end{align*}
where $v$ is the valuable population, $p$ is the pest population, $r$ is a population interaction constant which measures the efficiency of the pest population, $c$ is the intrinsic growth rate of the valuable population, $\kappa$ is the carrying capacity, $d$ is the intrinsic growth rate of the pest population, and $a$ is the control action to decrease the growth of the pest population. The parameters are chosen to be: $r = 0.5, c= 0.2, \kappa = 2, d= 0.2$. The objective is to steer the valuable population to $v = 1$ with the pest population and the control being $p = 0.2$ and $u=0.2$. The constraints are: $0.5 \le v\le 1.5$, $0\le p \le 1$, and $0\le u\le 0.4$. Let $x_1 = v-1,x_2 = p - 0.2, u = a-0.2$. Hence, $X = \{(x_1,x_2): \|x_1\|\le 0.5, -0.2 \le x_2 \le 0.8\}$ and $U=\{u: \|u\| \le0.2\}$. Let $N = 10,Q = I$ and $R = 0.1$. 

\begin{figure}[h]
\centering
\begin{minipage}[t]{.45\textwidth}
        \centering
		\includegraphics[width=0.9\linewidth]{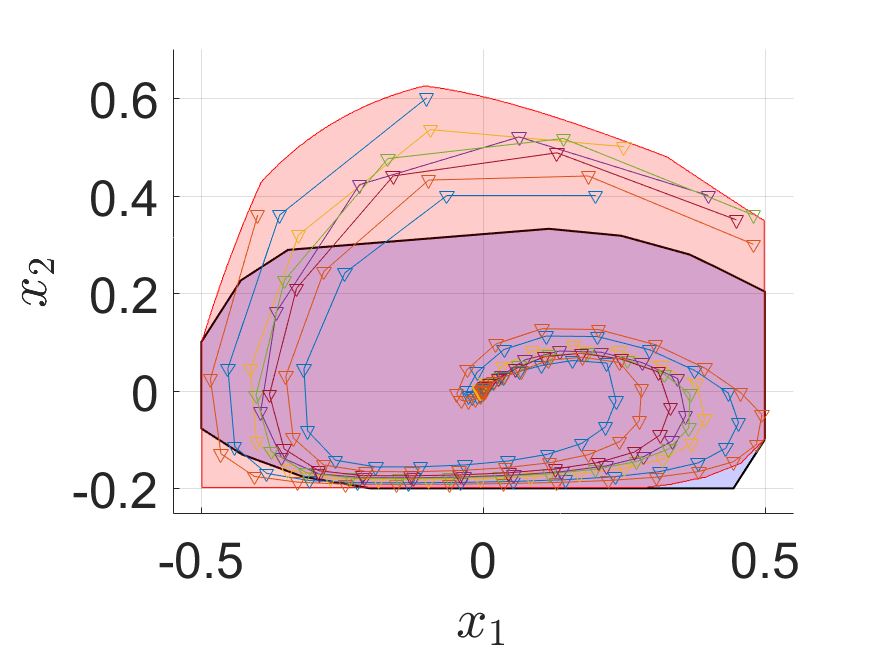}
		\caption{\emph{Polyflow} approximation versus Jacobian linearization of pest control dynamics: the blue area refers to $\mathcal{D}_{\Sigma(\bar{A},\bar{B},I)}$ and the red area refers to $\mathcal{D}_{\Sigma(A,B,C)}$.} 
\label{fig:domain}
\end{minipage}
\begin{minipage}[t]{0.45\textwidth}
\centering
		\includegraphics[width=0.9\linewidth]{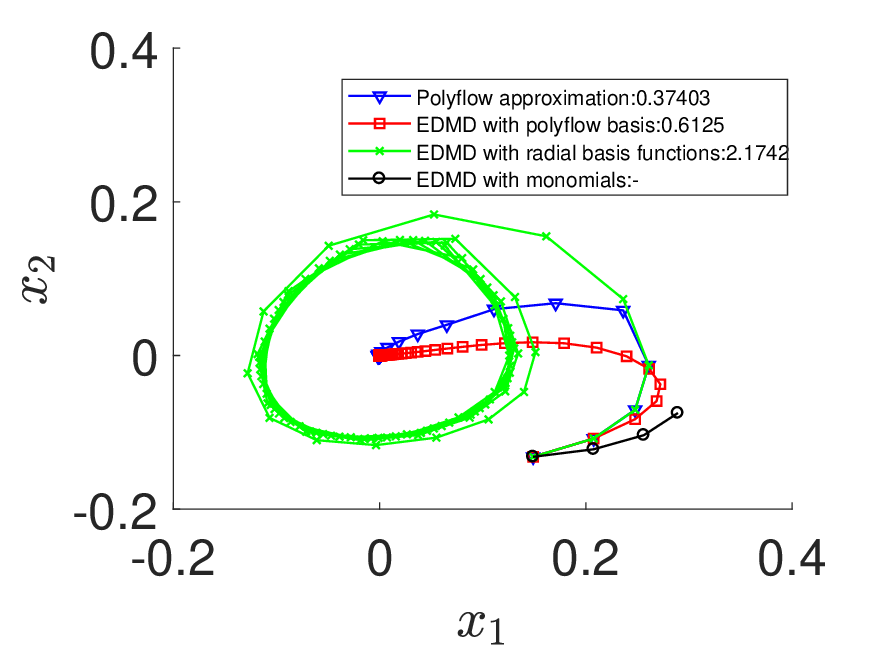}
		\caption{Comparison with radial and mononomial basis functions on pest control.} 
\label{fig:comp}
\end{minipage}
\end{figure}

First, we make comparison to Jacobian linearization about the origin $(x,u) = (0,0)$. Let $\bar{A} = \frac{\partial f(x,u)}{\partial x}|_{(0,0)}$ and $\bar{B} = \frac{\partial f(x,u)}{\partial u}|_{(0,0)}$. We then define $\tilde{\mathbb{P}}_{\Sigma(\bar{A},\bar{B},I)}$ as in (\ref{eqn:MPCABC}) with the feasible domain $\mathcal{D}_{\Sigma(\bar{A},\bar{B},I)}$. For the \emph{polyflow} approximation, we randomly sample $10^5$ points in $X$, solve Problem (\ref{eqn:alphakplus1}) with $k = 5$, and obtain a $5^{th}$ \emph{polyflow} approximation $\Sigma(A,B,C)$ with the feasible domain $\mathcal{D}_{\Sigma(A,B,C)}$. The simulation results are given in Figure \ref{fig:domain}. While recursive feasibility of the MPC problem is not guaranteed in both methods because of mismatch between the linearized model and the actual system , this figure shows that the \emph{polyflow} approximation outperforms Jacobian linearization as the \emph{polyflow} approximation is able to generate feasible trajectories from initial states that are not contained in $\mathcal{D}_{\Sigma(\bar{A},\bar{B},I)}$.

We also make comparison with Koopman-based approaches that use monomials \cite{ART:MG20} and radial basis functions.  We use monomials with the maximal degree being $6$, given by
$$
\{x_1^{s_1}x_2^{s_1}:1 \le s_1+s_2\le 6, s_1,s_2\in \mathbb{Z}^+\}.
$$
The number of monomials is $8!/(6!2!)-1=27$. The radial basis functions are taken from \cite{ART:KM18} in the form of
$$
g(x) = \|x-x_0\|^2\log(\|x-x_0\|),
$$
where $x_0$ is randomly selected with the uniform distribution on $X$. For a fair comparison, we generate $25$ radial basis functions and include the state itself in the basis. Hence, the dimension of the lifted state-space is also $27$. The matrices $A,B$ are obtained from Problem (\ref{eqn:ABM}). We then solve the MPC problems starting from the same initial state. To measure their performance quantitatively, we compute the LQ cost for $100$ steps:
$$
 \sum_{t=0}^{100} (\|x(t)\|_Q^2 + \|u(t)\|_R^2).
$$
The trajectories starting from $x(0) = [0.1488 ~ -0.1319]^T$ are shown in Figure \ref{fig:comp} with their corresponding LQ costs. As we can see from this figure, both the \emph{polyflow} approximation and the EDMD with \emph{polyflow} basis outperform the other two EDMD approaches. The LQ cost of the EDMD with monomials is not computed because it loses feasibility at $t=3$. The performance of the EDMD with radial basis functions relies on the randomly generated points $x_0$. Here, we select one realization where the MPC problem does not lose feasibility. Note that the dimension of the lifted system is $12$ for both the \emph{polyflow} approximation and  the EDMD with \emph{polyflow} basis.

\section{Conclusions}
In this paper, we present an approach to convert the non-convex MPC problem of nonlinear systems into a convex problem via linear embedding. The computation of linear embedding is implemented offline via \emph{polyflow} approximation, based on a nilpotency algebraic condition of the immersibility property for nonlinear control systems. From the derivation of this algebraic condition, we also show the link between the immersibility property and the Koopman operator. With this link, the proposed \emph{polyflow} approximation can be considered as a Koopman-based lifting approach with a special basis, which is computed endogenously. Finally, we evaluate our immersion-based MPC approach on the pest control problem by comparison with locally linearized MPC and other Koopman-based control approaches using monomial and radial basis functions.

\section*{Appendix}
\subsection*{Proof of Proposition  \ref{prop:TfABC}}
Sufficiency is obvious. We only need to prove necessity. Suppose System (\ref{eqn:fx}) is immersible into  $\Sigma(A,B, C)$. From Definition \ref{def:immersion}, there exists $T:\mathbb{R}^n \rightarrow \mathbb{R}^{\tilde{n}}$ such that, $\forall (x,u) \in \mathbb{R}^n \times \mathbb{R}^m$,
\begin{align*}
CT(x) &= x, \\
C(AT(x)+Bu) &= f(x,u) = CT(f(x,u)),\\
C(A^2T(x)+ABu+Bu) &= f(f(x,u),u) \\
& = C(AT(f(x,u))+Bu)
\end{align*}
Hence, $C(AT(x)+Bu-T(f(x,u))) = 0$ and $CA(AT(x)+Bu-T(f(x,u))) = 0$. Similarly, we can show that $CA^k(AT(x)+Bu-T(f(x,u))) = 0$ for all $k\in  \mathbb{Z}^+$. When $(C,A)$ is observable, we conclude that $AT(x)+Bu-T(f(x,u)) = 0$. If $(C,A)$ is not observable, we can find $(A',B',C')$ and $T'(x)$ in the observable subspace such that $A'T'(x)+B'u-T'(f(x,u)) = 0$ and $C'T'(x) = x$ hold following the same argument above. This completes the proof. $\Box$

\bibliographystyle{unsrt}
\bibliography{Reference}

\begin{thebibliography}{10}

\bibitem{ART:MRRS00}
D.~Q. Mayne, J.~B. Rawlings, C.~V. Rao, and P.~O.~M. Scokaert.
\newblock Constrained model predictive control: Stability and optimality.
\newblock {\em Automatica}, 36(6):789--814, 2000.

\bibitem{ART:C04}
M.~Cannon.
\newblock Efficient nonlinear model predictive control algorithms.
\newblock {\em Annual Reviews in Control}, 28(2):229--237, 2004.

\bibitem{ART:M14}
D.~Q. Mayne.
\newblock Model predictive control: Recent developments and future promise.
\newblock {\em Automatica}, 50(12):2967--2986, 2014.

\bibitem{ART:MN83}
S.~Monaco and D.~Normand-Cyrot.
\newblock The immersion under feedback of a multidimensional discrete-time
  non-linear system into a linear system.
\newblock {\em International Journal of Control}, 38(1):245--261, 1983.

\bibitem{ART:LM86}
J.~Levine and R.~Marino.
\newblock Nonlinear system immersion, observers and finite-dimensional filters.
\newblock {\em Systems \& Control Letters}, 7(2):133--142, 1986.

\bibitem{ART:LM88}
H.G. Lee and S.~I. Marcus.
\newblock Immersion and immersion by nonsingular feedback of a discrete-time
  nonlinear system into a linear system.
\newblock {\em IEEE transactions on Automatic Control}, 33(5):479--483, 1988.

\bibitem{ART:KI83}
A.~J. Krener and A.~Isidori.
\newblock Linearization by output injection and nonlinear observers.
\newblock {\em Systems \& Control Letters}, 3(1):47--52, 1983.

\bibitem{ART:AO03}
A.~Astolfi and R.~Ortega.
\newblock Immersion and invariance: A new tool for stabilization and adaptive
  control of nonlinear systems.
\newblock {\em IEEE Transactions on Automatic control}, 48(4):590--606, 2003.

\bibitem{ART:BM85}
H.~Bass and G.~Meisters.
\newblock Polynomial flows in the plane.
\newblock {\em Advances in Mathematics}, 55(2):173--208, 1985.

\bibitem{ART:VDE94}
A.~van~den Essen.
\newblock Locally finite and locally nilpotent derivations with applications to
  polynomial flows, morphisms and $\mathcal{G}_a$-actions. ii.
\newblock {\em Proceedings of the American Mathematical Society},
  121(3):667--678, 1994.

\bibitem{ART:VDE95}
A.~van~den Essen.
\newblock Locally nilpotent derivations and their applications, iii.
\newblock {\em Journal of Pure and Applied Algebra}, 98(1):15--23, 1995.

\bibitem{INP:JT19}
R.~M. Jungers and P.~Tabuada.
\newblock Non-local linearization of nonlinear differential equations via
  polyflows.
\newblock In {\em Proceedings of the American Control Conference}, 2019.

\bibitem{ART:WJ20}
Z.~Wang and R.~M. Jungers.
\newblock A data-driven immersion technique for linearization of discrete-time
  nonlinear systems.
\newblock {\em IFAC World Congress}, 2020.

\bibitem{ART:KKB17}
E.~Kaiser, J.~N. Kutz, and S.~L. Brunton.
\newblock Data-driven discovery of koopman eigenfunctions for control.
\newblock {\em arXiv preprint arXiv:1707.01146}, 2017.

\bibitem{ART:KM18}
M.~Korda and I.~Mezi{\'c}.
\newblock Linear predictors for nonlinear dynamical systems: Koopman operator
  meets model predictive control.
\newblock {\em Automatica}, 93:149--160, 2018.

\bibitem{ART:MG20}
A.~Mauroy and J.~Goncalves.
\newblock Koopman-based lifting techniques for nonlinear systems
  identification.
\newblock {\em IEEE Transactions on Automatic Control}, 65(6):2550--2565, 2020.

\bibitem{BOO:LVS12}
F.~L. Lewis, D.~Vrabie, and V.~L. Syrmos.
\newblock {\em Optimal control}.
\newblock John Wiley \& Sons, 2012.

\bibitem{ART:GT91}
E.G. Gilbert and K.~T. Tan.
\newblock Linear systems with state and control constraints: The theory and
  application of maximal output admissible sets.
\newblock {\em IEEE Transactions on Automatic Control}, 36:1008--1020, 1991.

\bibitem{INP:AJ16}
N.~Athanasopoulos and R.~M. Jungers.
\newblock Computing the domain of attraction of switching systems subject to
  non-convex constraints.
\newblock In {\em Proceedings of the 19th International Conference on Hybrid
  Systems: Computation and Control}, pages 41--50. ACM, 2016.

\bibitem{INP:WJO19}
Z.~Wang, R.~M. Jungers, and C.~J. Ong.
\newblock Computation of the maximal invariant set of linear systems with
  quasi-smooth nonlinear constraints.
\newblock In {\em Proceedings of the European Control Conference}, pages
  3803--3808, 2019.

\bibitem{ART:WJO21Com}
Z.~Wang, R.~M. Jungers, and C.~J. Ong.
\newblock Computation of the maximal invariant set of discrete-time linear
  systems subject to a class of non-convex constraints.
\newblock {\em Automatica}, 125:109463, 2021.

\bibitem{ART:WJO20}
Z.~Wang, R.~M. Jungers, and C.~J. Ong.
\newblock Computing invariant sets of discrete-time nonlinear systems via state
  immersion.
\newblock {\em IFAC World Congress}, 2020.

\bibitem{ART:WKR15}
M.~O. Williams, I.~G. Kevrekidis, and C.~W. Rowley.
\newblock A data--driven approximation of the koopman operator: Extending
  dynamic mode decomposition.
\newblock {\em Journal of Nonlinear Science}, 25(6):1307--1346, 2015.

\bibitem{INP:SM15}
I.~Susuki, Y.and~Mezi{\'c}.
\newblock A prony approximation of koopman mode decomposition.
\newblock In {\em Proceedings of the 54th IEEE Conference on Decision and
  Control}, pages 7022--7027, 2015.

\bibitem{ART:AM17}
H.~Arbabi and I.~Mezic.
\newblock Ergodic theory, dynamic mode decomposition, and computation of
  spectral properties of the koopman operator.
\newblock {\em SIAM Journal on Applied Dynamical Systems}, 16(4):2096--2126,
  2017.

\bibitem{ART:DHCLG14}
W.~Ding, R.~Hendon, B.~Cathey, E.~Lancaster, and R.~Germick.
\newblock Discrete time optimal control applied to pest control problems.
\newblock {\em Involve, a Journal of Mathematics}, 7(4):479--489, 2014.

\end{thebibliography}

\end{document}